# On The Equality Of The Grundy Numbers Of A Graph

Ali Mansouri and Mohamed Salim Bouhlel

Department of Electronic Technologies of Information and Telecommunications Sfax, Tunisia

*ABSTRACT*

*Our work becomes integrated into the general problem of the stability of the network ad hoc. Some, works attacked(affected) this problem. Among these works, we find the modelling of the network ad hoc in the form of a graph. Thus the problem of stability of the network ad hoc which corresponds to a problem of allocation of frequency amounts to a problem of allocation of colors in the vertex of graph. we present use a parameter of coloring " the number of Grundy". The Grundy number of a graph G, denoted by $\Gamma(G)$, is the largest k such that G has a greedy k-coloring, that is a coloring with colours obtained by applying the greedy algorithm according to some ordering of the vertices of G. In this paper, we study the Grundy number of the lexicographic, Cartesian and direct products of two graphs in terms of the Grundy numbers of these graphs.*

## KEYWORDS

*colouring, greedy algorithm, graph product ,grundy number, on-line algorithm,*

## 1. INTRODUCTION

Graphs considered in this paper are undirected, finite and contain neither loops nor multiple edges (unless stated otherwise). The definitions and notations used in this paper are standard and may be found in any textbook on graph theory. A (proper) k-colouring of a graph G=(V,E) is a mapping c:V →{1,...,k}, such that for any edge uv ∈ E(G),c(u) ≠c(v).

A k-colouring may also be seen as a partition of the vertex set of G into k disjoint stable set. $S_i$={v |c(v)=i} for 1≤i≤k. For convenience , by k-colouring we mean either the mapping c or the partition($S_1$,...,$S_k$).The elements of {1,...,k} are called colours. A graph is k-colourable if it has a k-colouring. The chromatic number χ(G) is the least k such that G is k-colourable. Several online algorithms producing colourings have been designed.

The most basic and most widespread one is the greedy algorithm. A greedy colouring relative to a vertex ordering σ=$v_1$<$v_2$<⋯<$v_n$ of V(G) is obtained by colouring the vertices in the order $v_1$,...,$v_n$, assigning to vi the smallest positive integer not already used on its lowered-indexed neighbours. Denoting by $S_i$ the stable set of vertices coloured i, a greedy colouring has the following property:

For every j<i, every vertex in Si has a neighbour in $S_j$ (1)

Otherwise the vertex in S would have been coloured j. Conversely, a colouring satisfying "Property (1)" is a greedy colouring relative to any vertex ordering in which the vertices of $S_i$ precede those of $S_j$ when i < j. The Grundy number $\Gamma(G)$ is the largest k such that G has a greedy k-colouring. It is well known that $\omega(G) \leq \chi(G) \leq \Gamma(G) \leq \Delta(G)+1$ where $\omega(G)$ denotes the clique number of G and $\Delta(G)$ the maximum degree of G.

The inequality $\chi(G) \leq \Gamma(G)$ may be tight, but it can also be very loose. Zaker [15] showed that for any fixed k≥ 0, given a graph G it is CoNP-Complete to decide whether $\Gamma(G) \leq \chi(G)+k$. He also showed that, given a graph G which is the complement of bipartite graph, it is CoNP-Complete to decide if $\Gamma(G)=\chi(G)$. This implies that it is CoNP-Complete to decide if $\Gamma(G)=\omega(G)$. Indeed, if G is the complement of a bipartite graph, then it is perfect, so $\chi(G)=\omega(G)$.

The Grundy number of various classes of graphs has been studied (see the introduction of [1]).

In this paper, we study the grundy number of different usual products of two graphs G and H. The lexicographic product G[H], the direct product G×H, and the Cartesian product G □ H, of G by H are the graphs with vertex set V(G) ×V(H) and the following edge set:

E(G[H]) = {(a,x) (b,y ) | ab ∈ E(G), or a=b and xy ∈ E(H)};

E(G×H) = {(a,x) (b,y) | ab ∈ E(G), and xy ∈ E (H)};

E(G□H) ={(a,x) (b,y) | a=b and xy ∈ E(H) or ab ∈ E(G) and x=y}.

It follows from the definition that G ×H (resp. G□H ) and H×G (resp. H□G ) are isomorphic. But G[H] and H[G] are generally not isomorphic. Moreover G[H] may be seen as the graph obtained by blowing up each vertex of G into a copy of H.

Regarding the lexicographic product, we prove in Section 3 that for any graphs G and H

$\Gamma(G) \times \Gamma(H) \leq \Gamma(G[H]) \leq 2^{\Gamma(G)-1}(\Gamma(H)-1)+\Gamma(G)-1$

In addition, we show that if G is a tree or $\Gamma(G)=\Delta(G)+1$, then $\Gamma(G[H])=\Gamma(G) \times \Gamma(H)$. Using these results, we prove a stronger complexity result than the one of Zaker [15] mentioned above: for every fixed c ≥1, it is CoNP-Complete to decide if $\Gamma(G) \leq c \times \chi(G)$ for a given graph G.

Analogously, we show that it is CoNP-Complete to decide if $\Gamma(G) \leq c \times \omega(G)$.

In Section 4, we investigate the Grundy number of the cartesian product of two graphs. We show that $(G \square H) \geq \max\{\Gamma(G), \Gamma(H)\}$ and increase this lower bound in some particular cases. We prove that there is no upper bound of $\Gamma(G \square H)$ as a function of $\Gamma(G)$ and $\Gamma(H)$. More precisely we show that for the complete bipartite $K_{p,p}$, $\Gamma(K_{p,p})=2$ but $\Gamma (K_{p,p} \square K_{p,p}) \geq p+1$.

Nevertheless, we show that for any fixed graph G, there is a function $h_G$ such that, for any graph H, $\Gamma(G \square H) \leq h_G(\Gamma(H))$.

in fact, we show that $h_G(k) \leq \Delta(G) \cdot 2^{k-1}+k$

We then give a better upper bound for $h_G(2)$ for some graphs G.

Finally, in Section 5, we study the Grundy number of the direct product of two graphs. We show that $\Gamma (G \times H) \geq \Gamma(G)+\Gamma(H)-2$ and construct for any k some graph $G_k$ such that $\Gamma(G_k)=2k+1$ and $\Gamma(G_k \times K_2)=3k+1$.

## 2. PRELIMINAIRES

In this section, we present some definitions and preliminary results.

A sub graph of a graph G is a graph H such that V (H) ⊂V (G) and E(H)⊂E(G). Note that since H is a graph we have E (H) ⊂ E (G) ∩ [V (H)] $_2$. If H contains all the edges of G between vertices of V (H), that is E (H) =E (G) ∩V (H)] $_2$, then H is the subgraph induced by V (H).

If S is a set of vertices, we denote by G<S> the graph induced by S and by G−S the graph induced by V (G)\S.

For simplicity, we write G −v rather than G−{v}. For a subset F of E(G), we write G\F
= (V (G), E (G)\F).

For a subset F of E (G), we write G\F = (V (G), E (G)\F). As above G \{e}is abbreviated to G\e.

If H is a subgraph of G then χ(H)≤χ(G). This assertion cannot be transposed to the grundy number. For example, the path $P_4$ of order 4 is a subgraph of the cycle $C_4$ of order 4 but one can easily check that Γ($P_4$)=3 and Γ($C_4$)=2. However such an assertion holds if we add the extra condition of being an induced subgraph.

**Proposition1** If H is an induced subgraph of G then Γ(H)≤Γ(G).

**Proof**.

Let σ be an ordering for which the corresponding greedy colouring of H uses Γ (H) colours. Then a colouring with respect to any ordering of V (G) beginning with σ will use at least Γ (H) to colour H, hence at least Γ (H) to colour G.

**Lemma2** Let G be a graph and u and v two vertices G. The following hold:

For any edge e, Γ(G) −1≤Γ(G\e)≤Γ(G)+1. (1)

If N (u) ⊂ N (v) then in every greedy colouring c of G, c(u)≤c(v). In particular, if N (u) =N (v) then:

c (u)=c(v). (2)

If N (u) =N (v) then Γ (G) =Γ (G−u). (3)

**Proof.**

(1) Set e= xy and p= Γ (G). Let ($S_1$,...,$S_p$) be a greedy p-colouring of G. It satisfies Property (1).

Let i be the integer such that x ∈ $S_i$ and let $T_j$=$S_j$ for 1≤ j< i and $T_j$=S $_{j+1}$ for i≤ j≤ p−1.

It is a simple matter to check that ($T_1$,...,T $_{p-1}$) satisfies Property (1). Hence Γ (G−$S_i$) ≥ p-1. As G−$S_i$ is an induced subgraph of G\e, by proposition1, Γ (G\e) ≥ p−1.

Set q=Γ (G). Let ($S_{1'}$,...,$S_{q'}$) be a greedy q-colouring of G\e. It satisfies Property (1). Now let i be the integer such that x ∈ $S_{i'}$. Let $T_{j'}$= $S_{j'}$ for 1≤ j<i and $T_{j'}$= S $_{j+1'}$ for i≤ j≤q−1. It is a simple matter to check that ($T_{1'}$,..., T $_{q-1'}$) satisfies Property (1). Hence Γ(G−$S_i$)≥ q−1. As G−$S_i$ is an induced subgraph of G, by Proposition 1, Γ (G) ≥q−1.

(2) Let c= ($S_1$,...,$S_p$). Suppose u ∈$S_j$ and v ∈ S $_i$. Since v ∈ $S_i$, then v has no neighbour in $S_i$. So u has no neighbour in $S_j$ because N (u) ⊂N (v). Thus j≤ i because c satisfies Property (1).

(3) Let $S_1$,...,$S_p$ be the stable sets of a greedy coulouring..

By (2), u and v are in the same stable set $S_i$.

Now $S_1$,...,S $_{i−1}$,$S_i$\{u},S $_{i+1}$,...,$S_p$ are the stable sets of a greedy colouring of G−u. Indeed as $N_G$ (u) =$N_G$ (v) it is a simple matter to check that they satisfy Property (1).

A path is a non-empty graph P= (V, E) of the form V= {$x_0, x_1, ..., x_k$} and E={$x_0x_1, x_1x_2, ..., x_{k-1}x_k$} where the $x_i$ are all distinct. The vertices $x_0$ and $x_k$ are the endvertices of P. A (u,v)-path is a path with endvertices u and v. A graph is connected if for any two vertices u an v there is a (u,v)- path.

**Proposition3** Let G be a connected graph. Then $\Gamma(G)=2$ if and only if G is complete bipartite.

**Proof**. It is easy to see that if G is complete bipartite then $\Gamma(G)=2$: indeed applying several times.

**Lemma 2** (3), we obtain that $\Gamma(G) = \Gamma(K_2) = 2$.

Conversely, if $\Gamma(G) = 2$, then G has to be bipartite because $\Gamma(G) \geq \chi(G)$. Suppose now that G is not complete bipartite. Then there exist two vertices u and v in different parts of the partition which are not adjacent. Let P be a shortest (u,v)-path. Then P has odd length, so length at least 3 and because it is a shortest path it is an induced path. Hence G contains an induced $P_4$. So by Proposition 1, $\Gamma(G) \geq 3$.

This proposition implies that one can decides in polynomial time if the Grundy number of a graph is 2. More generally, Zaker [15] showed that for any fixed k, it is decidable in polynomial time if a given graph has grundy number at most k. To show this, he proved that there is a finite number of graphs called k-atoms such that if $\Gamma(G) \geq k$ then G contains a k-atom as induced subgraph. The k-atoms may easily be found using Proposition 5 below.

**Definition4** Let G be a graph and W a subset of V (G). A set S is W-dominating if $S \subset V(G) \setminus W$ and every vertex of W has a neighbour in S.

The following proposition follows immediately from the Property (1) of greedy colouring.

**Proposition5** Let G be a graph and W a subset of V (G). If S is a W-dominating stable set then $\Gamma(G<W \cup S>) \geq \Gamma(G<W>)+1$.

Note that if S is a W-dominating set then $\Gamma(G<W \cup S>)$ cannot be bounded by a function of $\Gamma(G<W>)$.

For example, a tree may be partitioned into two stable sets S and T. Moreover, because the tree is connected S is T-dominating (and vice-versa). But the Grundy number of a stable set is 1 whereas the Grundy number of a tree may be arbitrarily large. Consider for example the binomial tree of index k $T_k$ which may be defined recursively as follows:

- $T_1$ is the graph with one vertex and no edge;
- $T_k$ is constructed from $T_{k-1}$ by joining each vertex to a new leaf.

The binomial tree $T_k$ has chromatic number 2 and grundy number k. It is the unique k-atom which is a tree. Hence, as shown in [8], the Grundy number of a tree is the largest index of a binomial tree it contains.

The union of two graphs $G_1$ and $G_2$ is the graph $G_1 \cup G_2$ with vertex set $V(G_1) \cup V(G_2)$ and edge set $E(G_1) \cup E(G_2)$. If $G_1$ and $G_2$ are disjoint (i.e. $V(G_1) \cap V(G_2) = \phi$), we refer to their union as a disjoint union and denote it $G_1+G_2$. The join of two disjoint graphs $G_1$ and $G_2$ is the graph $G_1 \oplus G_2$ obtained from $G_1+G_2$ by joining all the vertices of $G_1$ to all the vertices of $G_2$.

**Proposition6**   If $G=G_1+G_2$ then $\Gamma(G) = \max(\Gamma(G_1), \Gamma(G_2))$.
   If $G = G_1 \oplus G_2$ then $\Gamma(G) = \Gamma(G_1) + \Gamma(G_2)$.

This proposition and an immediate induction yield a result of Gyarfás and Lehel [6] stating that for every cograph (graph without induced $P_4$) $\Gamma(G) = \chi(G)$ because every cograph of order at least two is either the disjoint union or the join of two cographs.

**Lemma7** Let G be a graph and x a vertex of G. If there is a greedy colouring c such that x is coloured p then for any $1 \leq i \leq p$, there is a greedy colouring such that x is coloured i.

**Proof**. For $1 \leq i \leq p-1$, let $S_i$ be the stable set of vertices coloured i by c. Then for any $1 \leq i \leq p$, $(S_1,...,S_{i-1},\{x\})$ is a greedy i-colouring of $G\langle\{x\}\cup\bigcup_{j=1}^{i-1} S_i\rangle$ in which x is coloured i. This partial greedy colouring of G may be extended into a greedy colouring of G in which x is coloured i.

**Lemma8** Let G be a graph with at least one edge. There are two adjacent vertices x and y such that there is two greedy colourings $c_x$ and $c_y$ such that $c_x(x) = c_y(y) = \Gamma(G)$.

**Proof**. Set $p = \Gamma(G)$ and let $c_x$ be a greedy p-colouring of G with stable sets $S_1,...,S_p$. Let x a vertex of $S_p$ and y a neighbor of x in $S_{p-1}$. Then $S_1,S_2,...S_{p-2},S_p,\{y\}$ is a partial greedy colouring $c_y$ of G with $c_y(x)= p-1$ and $c_y(y)= p$. This colouring may trivially be extended to G.

## 3. LEXICOGRAPHIC PRODUCT

Obviously, $\chi(G[H]) \leq \chi(G) \times \chi(H)$ and Stahl [11] showed $\chi(G[H]) \geq \chi(G)+2\chi(H)-2$. In this section, we establish some bounds on $\Gamma(G[H])$ in terms of $\Gamma(G)$ and $\Gamma(H)$.

**Definition9** [12] In the lexicographic product G[H], for every vertex $x \in G$, we call copy of H at x the graph H(x) isomorphic to H which is induced by the vertices of $\{x\} \times V(H)$.

**Proposition10** Let G and H be two graphs. In a greedy colouring of G[H], at most $\Gamma(H)$ colours appear on each H(x), $x \in V(G)$.

**Proof**. Consider a greedy colouring of G[H] and let $n_1,n_2,...,n_p$ be the p colours appearing on a particular copy H(x) of H. For any $1 \leq i \leq p$, let $S_i$ be the stable set of vertices of H(x) coloured $n_i$. Let u be a vertex of $S_i$. For any $1 \leq j < i$, by the Property (1), in G[H], u has a neighbour v coloured $n_j$. The vertex v must be in H(x) because the neighbours of x not in H(x) are also neighbours of the vertex z of H(x) coloured $n_j$. Hence $v \in S_j$. It follows that the colouring $(S_1,...,S_p)$ satisfies the Property (1). Hence $\Gamma(H)= \Gamma(H(x)) \geq p$. Geller and Stahl [5] showed that if $\chi(H)=k$ then $\chi(G[H])=\chi(G[K_k])$ for any graph G. We now prove a similar result for the Grundy number.

**Theorem11** Let H be a graph such that $\Gamma(H)=k$. Then for any graph G, $\Gamma(G[H]) = \Gamma(G[K_k])$.

**Proof**. Set $V(G) = \{v_1,...,v_n\}$.

Let c be a greedy colouring of G[H]. For every $1 \leq i \leq n$, let $A_i = c(H(V_i))=\{\alpha^1_i,...,\alpha^{|A_i|}_i\}$ be the set of colours appearing on $H(v_i)$. Let F be the graph obtained from G[H] by replacing each $H(v_i)$ by a complete graph on $|A_i|$ vertices, $w^1_i,...,w^{|A_i|}_i$ and c' be the colouring of F defined by $c'(w^j_i)= \alpha^j_i$ for any $1 \leq i \leq n$ and $1 \leq j \leq |A_i|$.

By construction F is an induced subgraph of $G[K_k]$ because for each $1 \leq i \leq n$, $|A_i| \leq k$ by Proposition10. Moreover, it is a simple matter to check that c' is a greedy colouring of F. Hence $\Gamma(G[K_k]) \geq \Gamma(F) \geq \Gamma(G[H])$.

### 3.1. Lower bounds

**Proposition 12** Let G and H be two graphs. Then $\Gamma(G[H]) \geq \Gamma(G) \times \Gamma(H)$.

**Proof of Proposition 12**. Let $c_G$ (resp. $c_H$) be a greedy colouring of G (resp. H) with $\Gamma(H)$ (resp. $\Gamma(G)$) colours. Then the colouring $c=(c_G,c_H)$ with the pairs of colours ordered according to the lexicographic product is a greedy colouring of G[H].

International Journal of Next-Generation Networks (IJNGN) Vol.5, No.4, December 2013

Proposition 12 is tight as there are pairs of graphs (G, H) for which $\Gamma(G[H]) = \Gamma(G) \times \Gamma(H)$. In particular, we shall prove that if G is a tree or satisfies $\Gamma(G)=\Delta(G)+1$ this is the case.

**Theorem 13** Let G be and H be two graphs. If $\Gamma(G)=\Delta(G)+1$ then $\Gamma(G[H])=\Gamma(G)\times\Gamma(H)$.

**Proof.** By Proposition 12, $\Gamma(G[H]) \geq \Gamma(G) \times \Gamma(H)$.

Let us now show that $\Gamma(G[H]) \leq \Gamma(G) \times \Gamma(H)$. Consider a greedy colouring of G[H]. Let u be a vertex of G[H] coloured with the largest colour $c_{max}$ and H(x) the copy of H containing u. Moreover by the the Property (1), every colour but $c_{max}$ must appear on the neighbour hood of u. Hence $c_{max} \leq \Gamma(G) \times \Gamma(H)$.

**Theorem 14** Let T be a tree and H be an arbitrary graph. Then $\Gamma(T[H])=\Gamma(T)\times\Gamma(H)$.

**Proof.** Let k be the integer such that $k\Gamma(H) \geq \Gamma(T[H]) \geq (k-1)\Gamma(H)+1$. We will prove that $\Gamma(T) \geq k$ by showing that T contains a binomial tree of index k as an induced subgraph. This implies that $\Gamma(T[H]) \leq \Gamma(T) \times \Gamma(H)$.

So by Proposition 12, $\Gamma(T[H]) = \Gamma(T) \times \Gamma(H)$.

Let f be a greedy colouring of T[H] with $\Gamma(T[H])$ colours in the following, by colour we should understand colour assigned by f. We shall construct step by step a binomial tree of order k in T.

**Step 1**: Let $v_1$ be a vertex of T such that a vertex of $H(v_1)$ is coloured $c_1=\Gamma(T[H])$. Then the sub tree of T with unique vertex $v_1$ is $T_1$. Let $P_1(v_1)$ be the sequence $(v_1)$.

Step i: (for $2 \leq i \leq k$) we have the binomial sub tree $T_{i-1}$ of T. Moreover, to each vertex v of $T_{i-1}$ is associated a sequence $P_{i-1}(v)=(v_{i-1}, v_{i-2},...,v_2,v_1)$ of $i-1$ vertices in $T_{i-1}$ such that :

(a) $P_{i-1}(v)$ contains v and all its neighbours in $T_{i-1}$, and
(b) for all $1 \leq j \leq i-1$, $H(v_j)$ contains the greatest colour not appearing on $\cup_{l=1}^{j-1} H(v_l)$.

We shall construct $T_i$, that is add a leaf to each vertex of $T_{i-1}$, and also describe the sequences $P_i$ satisfying the conditions (a) and (b). Let v be a vertex of $T_{i-1}$.

As $P_{i-1}(v)$ contains $i-1$ vertices, at most $(i-1)\Gamma(H)$ colours appear on $\cup_{l=1}^{i-1} H(v_l)$ by Proposition 10.

### 3.2. Upper bounds

There are pairs of graphs (G, H) for which $\Gamma(G[H])$ is greater than $\Gamma(G) \times \Gamma(H)$ as shown by the following proposition.

**Proposition 15** Let $G_3$ be the graph depicted. Then $\Gamma(G_3)=3$ and $\Gamma(G_3[K_{2p}]) \geq 7p$

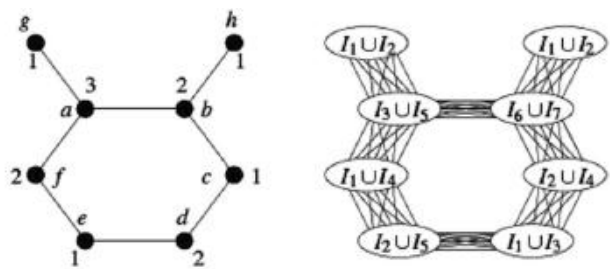

Figure 1: A greedy 3-colouring of $G_3$ and a greedy 7p-colouring of $G_3[K_{2p}]$.

**Proof**. Let us first show that $\Gamma(G_3) = 3$.

Suppose, by way of contradiction, that $G_3$ admits a greedy 4-colouring. Then one of the two vertices of degree three, namely a and b, is coloured 4. By symmetry, we may assume that it is a. This vertex must have a neighbour coloured 3. This neighour is necessarily b which is the unique one having degree at least two in $G_3 - a$. The vertices a and b must each have a neighbour coloured 2 which must have degree at least one in $G - \{a,b\}$. Hence f and c are coloured 2. These two vertices must have a neighbour coloured 1. So d and e are coloured 1, which is a contradiction as they are adjacent.
[

Let us now show that $\Gamma(G_3[K_{2p}]) \geq 7p$. For every vertex $v \in V(G_3)$, let us assign 2p colours to the 2p vertices of $K_2(v)$ as follows. $I_3 \cup I_5$ to $K_2(a)$, $I_6 \cup I_7$ to $K_2(b)$, $I_2 \cup I_4$ to $K_2(c)$, $I_1 \cup I_3$ to $K_2(d)$, $I_2 \cup I_5$ to $K_2(e)$, $I_1 \cup I_4$ to $K_2(f)$, $I_1 \cup I_2$ to $K_2(g)$ and $K_2(h)$. It is a simple matter to check that this is a greedy 7p-colouring of $G_3[K_{2p}]$. We would like to find upper bounds on $\Gamma(G[H])$ in terms of
$\Gamma(G)$ and $\Gamma(H)$.

Ideally we would like to determine exactly
$\Psi(k,l) = \max\{\Gamma(G[H]) | \Gamma(G)=k \text{ and } \Gamma(H)=l\}$
$= \max\{\Gamma(G[K_l]) | \Gamma(G) = k\}$

by Theorem 11. In the remainder of this section we give upper and lower bounds on $\psi$. Note that $\Gamma(G) = 1$ if and only if G has no edge. Thus if $\Gamma(H) = 1$ then $\Gamma(G[H]) = \Gamma(G)$ using Lemma 2-(3). Moreover if $\Gamma(G) = 1$ then $\Gamma(G[H]) = \Gamma(H)$ by Proposition 6. In the remainder of the section, we will assume that all the graphs we consider have Grundy number at least 2.

**Theorem 16** $\Gamma(G[H]) \leq 2^{\Gamma(G)-1}(\Gamma(H)-1) + \Gamma(G) - 1$

**Proof**. Let k be the integer such that $2^{k-1}(\Gamma(H)-1) + k - 1 \geq \Gamma(G[H]) > 2^{k-2}(\Gamma(H)-1) + k - 2$.

We will show that $\Gamma(G) \geq k$, which implies that $\Gamma(G[H]) \leq 2^{\Gamma(G)-1}(\Gamma(H)-1) + \Gamma(G) - 1$.

Let f be a greedy colouring of $G[H]$ with $\Gamma(G[H])$ colours. We shall construct step by step an induced subgraph of G which has Grundy number at least k.

**Step 1**: Let $v_1$ be a vertex such that the largest colour $c_1 = \Gamma(G[H])$ appears on $H(v_1)$. Let $G_1 = G<\{v_1\}>$. Then $\Gamma(G_1) = 1$.

**Step 2**: Since $\Gamma(G[H]) > 2^{k-2}(\Gamma(H)-1) + k - 2 \geq \Gamma(H)$, by Proposition 10, there are colours that do not appear on $H(v_1)$. Let $c_2$ be the largest such colour. For $c_1 > c_2$, there is a vertex $v_2 \in N_G(v_1)$ such that $c_2$ appears on $H(v_2)$. Let $G_2 = G<\{v_1,v_2\}>$. Since $v_1v_2$ in an edge, $\Gamma(G_2) = 2$

**Step i**: (for $3 \leq i \leq k$): We have a subgraph $G_{i-1}$ of G of at most $2^{i-2}$ vertices such that $\Gamma(G_{i-1}) \geq i-1$ and at most $2^{i-2}(\Gamma(H)-1) + i - 2$ colours appear on $G_{i-1}[H]$. For $\Gamma(G[H]) > 2^{k-2}(\Gamma(H)-1) + k - 2 \geq 2^{i-2}(\Gamma(H)-1) + i - 2$, there are colours that do not appear on $G_{i-1}[H]$. Let $c_i$ be the greatest such colour. Since $c_0 > c_1 > ... > c_i$ and $c_i$ does not appear on $G_{i-1}[H]$, any vertex $v \in V(G_i)$ has a neighbour $n(v)$ in $V(G)\setminus V(G_i)$ such that the colour $c_i$ appears on $H(n(v))$. Let $S_i = \{n(v), v \in V(G_i)\}$ and $G_i = G<V(G_{i-1}) \cup S_i>$. Then $|S_i| \leq |G_{i-1}|$ so $|G_i| \leq 2^{i-1}$. Moreover $S_i$ is a stable set since the colour $c_i$ appears on the copy of H at each vertex of $S_i$. So by Proposition 5, $\Gamma(G_i) \geq \Gamma(G_{i-1}) + 1 \geq i$. Now at most $2^{i-2}(\Gamma(H)-1) + i - 2$ colours appear on $G_{i-1}[H]$ and at most $2^{i-2}(\Gamma(H)-1) + 1$ colours appear on $S_i[H]$ by Proposition 10 and because $c_i$ appears in all the $H(v)$ for $v \in S_i$. So in total at most $2^{i-1}(\Gamma(H)-1) + i - 1$ colours appear on $G_i[H]$.

**Corollary 17**

$\psi(k,l) \leq 2^{k-1}(l-1) + k - 1.$     (1)

If $\Gamma(G) = 2$ then $\Gamma(G[H]) = 2k.$     (2)

$\psi(2, k) = 2k.$ (3)

$(\psi(3,2) = 7.$ (4)

**Proof**. (1) follows directly Theorem 16; Proposition 12 and Theorem 16 imply (2) and (3); Proposition 15 and Theorem 16 yield (4).

**Lemma 18** Let $\alpha$ be a positive integer. If $\psi(k,l) \geq kl + \alpha$ then $\psi(k',l) \geq k'l + \alpha$ for all $k' > k$.

**Proof**. To prove this result it suffices to prove that if $\psi(k,l) = kl + \alpha$ then $\psi(k+1,l) \geq (k+1)l + \alpha$. Then an easy induction will give the result.

Let G be a graph such that $\Gamma(G[K_l]) = k_l + \alpha$. Let x be a vertex of G such that there exists a greedy $(kl + \alpha)$-colouring c such that the colour $kl + \alpha$ appears on $K_l(x)$. Let $G_1$ and $G_2$ be two disjoint copies of G. For $i = 1,2$, we denote by $v_i$ the vertex $v_i \in V(G_i)$ corresponding to $v \in V(G)$. Let G' be the graph obtained from $G_1+G_2$ by adding an edge between the two vertices $x_1$ and $x_2$. By Lemma 2 (1) and Proposition 6, $\Gamma(G') \leq \Gamma(G_1+G_2)+1 = \Gamma(G)+1 = k+1$. Now let c' be the colouring of $G'[K_l]$ defined as follows:

- $c'(v_1) = c(v)$ for $v_1 \in V(G_1[K_l])$;
- $c'(v_2) = c(v)$ for $v_2 \in V(G_2[K_l]) \setminus \{x_2\}$;
- The vertices of $K_l(x_2)$ are assigned distinct colours in $\{k_l + \alpha + 1, ..., (k+1)l + \alpha\}$.

**Corollary 19** Let $k \geq 3$ be an integer. Then $\psi(k,p) \geq (2k+1)p$

**Theorem 20** Let $c \geq 1$ be an integer. The following problem is CoNP-complete:
Instance: a graph G
Question: $\Gamma(G) \leq c\chi(G)$?

**Proof**. Let G be a graph. If $c_1$ is a colouring of G with t colours and $c_2$ a greedy colouring of G with more than ct colours, then the pair $(c_1,c_2)$ forms a certificate that $\Gamma(G) > c\chi(G)$. Clearly, it can be checked in polynomial time if a pair $(c_1,c_2)$ is a certificate. So the problem is in CoNP.

Let us now show that this problem is CoNP-complete via a reduction to the problem of deciding if $\Gamma(G) \leq \chi(G)$ for a given graph G, which is known to be CoNP-complete [15]. Let G be a graph.

Consider $H = T_{2c}[G]$. Then $\chi(H) = 2\chi(G)$ as $\omega(T_{2c}) = \chi(T_{2c}) = 2$. Moreover $\Gamma(H) = 2c\Gamma(G)$ by Theorem 13 (or Theorem 14). Hence $\Gamma(H) \leq c\chi(H)$ if and only if $\Gamma(G) \leq \chi(G)$.

A similar proof yields that it is NP-complete to decide if $\Gamma(G) \leq c\omega(G)$ as $\omega(T_{2c}[G]) = 2\omega(G)$ and it is CoNP-complete to decide if $\Gamma(G) \leq \omega(G)$.

**Theorem 21** Let $c \geq 1$ be an integer. The following problem is CoNP-complete:
• Instance : a graph G.
• Question : $\Gamma(G) \leq c\,\omega(G)$?

## 4. CARTESIAN PRODUCT

It is well-known that the chromatic number of the cartesian product of two graphs is the maximum of the chromatic numbers of these graphs: $\chi(G \square H) = \max\{\chi(G), \chi(H)\}$. Unfortunately, no such formula holds for the grundy number. In this section, we are looking for bounds on the grundy number of the cartesian product of two graphs in terms of the grundy numbers of these graphs. We first show that such an upper bound does not exist. However, we show that for any graph G there is a function $h_G$ such that for every graph H, $\Gamma(G \square H) \leq h_G(\Gamma(H))$. Regarding lower bounds, we give upper an lower bounds for the function

$\phi_\Box(k,l) = \min\{\Gamma(G\Box H)|\Gamma(G)=k \text{ and } \Gamma(H)=l\}$

Let G and H be two graphs. For any $v \in V(G)$, the graph $H_v$ of $G\Box H$ induced by the vertices of $\{v\}\times V(H)$ is isomorphic to H. Analogoulsy, for any $x \in V(H)$, the subgraph $G_x$ of $G\Box H$ induced by the vertices of $V(G)\times\{x\}$ is isomorphic to G.

### 4.1. Upper bounds

We denote by $K_{p,p}$ the complete bipartite graph with p vertices in each part.

**Proposition 22** Let p≥2 be an integer. Then $\Gamma(K_{p+1,p+1}\Box K_{p+1,p+1}) \geq \Gamma(K_{p,p}\Box K_{p,p})+1$.
So $\Gamma(K_{p,p}\Box K_{p,p}) \geq p+1$.

**Proof.** Let $(X\cup\{x\}, Y\cup\{y\})$ be the bipartition of $K_{p+1,p+1}$ with $x \not\in X$ and $y \not\in Y$.
Then $K_{p+1,p+1}-\{x,y\}$ is a $K_{p,p}$, so $K_{p+1,p+1}-\{x,y\}\Box K_{p+1,p+1}-\{x,y\}$ is an induced $K_{p,p}\Box K_{p,p}$ in $K_{p+1,p+1}\Box K_{p+1,p+1}$.

Now the set $(\{x\}\times Y\setminus\{y\}) \cup (\{y\}\times X\setminus\{x\}) \cup (X\setminus\{x\}\times\{x\}) \cup (Y\setminus\{y\}\times\{y\})$ is a $(X\cup Y)\times(X\cup Y)$ dominating stable set. So by Proposition 5, $\Gamma(K_{p+1,p+1}\Box K_{p+1,p+1}) \geq \Gamma(K_{p,p}\Box K_{p,p})+1$.

As $\Gamma(K_{2,2})=2$, an easy induction yields $\Gamma(K_{p,p}\Box K_{p,p}) \geq p+1$, for p≥2.

**Proposition 23** Let G be a graph then for any positive integer k, $h_G(k) \leq \Delta(G)\cdot 2^{k-1}+k$. In other words, for any graph H, $\Gamma(G\Box H) \leq \Delta(G)\cdot 2^{\Gamma(H)-1}+\Gamma(H)$.

**Proof.** Let c be a greedy p-colouring of $G\Box H$. Let $(v,x_1)$ be a vertex coloured $p=c_1$. For every vertex x of H, set $C(x):=\{c(w,x)|w \in N_G(v)\}$. By extension, for every $S\in V(H)$, we set $C(S) =\cup_{x\in S} C(x)$. Let $T_1=\{x_1\}$.

We have $\Gamma(H<T_1>)=1$. Now, iteratively, as long as $\{1,...,p\}\setminus C(T_i) \cup\{c_1,...,c_i\}$ is not empty, let us construct $T_{i+1}$ as follows. Let $c_{i+1}$ be the largest integer of $\{1,...,p\}\setminus C(T_i) \cup\{c_1,...,c_i\}$. Then for every $x \in T_i$, the vertex $(v,x)$ has a neighbour coloured $c_{i+1}$ which by definition of C(x) is in $H_v$. Hence there exists a stable set $S_{i+1}$ of size at most $|T_i|$ in H such that $c(v,y)=c_{i+1}$ for every $y \in S_{i+1}$ and every vertex $x \in T_i$ has a neighbour in $S_{i+1}$. Setting $T_{i+1}=T_i\cup S_{i+1}$, we have $|T_{i+1}|\leq 2|T_i|\leq 2^i$ and by Proposition 5, $\Gamma(H<T_{i+1}>) \geq i+1$

Let $i_0$ be the integer when the process terminates, i.e. when $\{1,...,p\}=C(T_{i0}) \cup\{c_1,...,c_{i0}\}$.
We have $\Gamma(H) \geq \Gamma(H<T_{i0}>) \geq i_0, |T_{i0}|\leq 2^{i_0-1}$ and $|C(T_{i0})|\leq \Delta(G)\times|T_{i0}|$. So $p\leq \Delta(G)\cdot 2^{i_0-1}+i_0 \leq \Delta(G)\cdot 2^{\Gamma(H)-1}+\Gamma(H)$.

We think that the upper bound $\Delta(G) \cdot 2^{k-1}+k$ is far to be tight. For some graphs one can get slightly better upper bounds.

Let us show an example when k=2. For a vertex v of graph G, we denote by $d_1G(v)$ or simply $d_1(v)$ the maximum degree of a neighbour of v, i.e. $d_1(v)=\max\{d(u)|u \in N(v)\}$.

According to the proof of Theorem 25, $p\leq \max\{d_G(v)+d_1G(v)+2 \mid v \in V(G)\}$. We now show a slightly better upper bound.

**Proposition 26** Let G be a graph. Then $h_G(2)\leq\max\{\min\{2d(v)+2, 2d_1(v)+3\}|v \in V(G)\}$.

**Proof.** Let H be a complete bipartite graph and c be a greedy colouring of $G\Box H$ with p colours. Let $x=(v, v')$ be a vertex coloured with p and let (X,Y) be the bipartition of $H_v$ with $x \in X$.

Since x has $d_G(v)$ neighbours not in $H_v$, it has $p-1-d_G(v)$ neighbours in Y with distinct colours in $\{1,...,p-1\}$. Let q be the largest integer in $\{1,...,p-1\}$ that is assigned to a vertex in Y and let y be a vertex coloured q. Then x has $p-2-d_G(v)$ neighbours in Y with distinct colours in $\{1,...,q-1\}$.

Now since y has at most $d_G(v)$ neighbours not in $H_v$, it has $q-1-d_G(v)$ neighbours in X with distinct colours in $\{1,...,q-1\}$. As $H_v$ is complete bipartite, the colours that appear on X do not appear on Y.

Thus $p-2-d_G(v)+q-1-d_G(v) \leq q-1$, so $p \leq 2d_G(v)+2$.

We claim that there is a vertex $y=(u, u')$ with $u \in N_G(v)$ such that is assigned a colour $p' \geq p-2$ and is adjacent to a vertex in $H_v$ coloured p or $p-1$. Indeed x has a neighbour that is coloured $p-1$.

If this neighbour is not in $H_v$ it is the desired y. If not this neighbour z is in Y. Now both x and z have a neighbour coloured $p-2$. But these two neighbours are not both in $H_v$ otherwise they would be adjacent. Hence one of them is not in $H_v$ and is the desired y.

Now applying the same reasoning as above and taking into account that y has a neighbor outside $H_u$ with a larger colour than it's, we obtain that $p-2 \leq 2d_G(u)+1$. So $p \leq 2d_1(v)+3$.

If the graph G has two adjacent vertices of G maximum degree then Proposition 26 yields the same upper bound $2\Delta(G)+2$ as Theorem 25. But for graphs in which vertices of high degree form a stable set, this bound is far better. Consider for example a star $K_{1,p}$. By Proposition 26, for any $p \geq 2$, $hK_{1,p}(2) \leq 5$. Moreover $K_{1,p}$ contains $K_{1,2}$ as an induced subgraph, so $K_{1,p} \square K_{3,3}$ contains $K_{1,2} \square K_{3,3}$ as an induced subgraph.

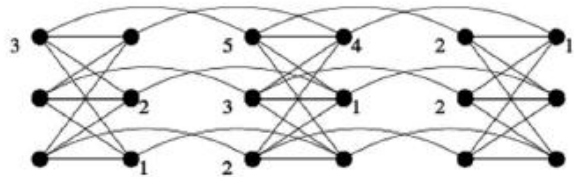

Figure 2: Partial greedy 5-colouring of $K_{1,2} \square K_{3,3}$

With similar arguments, one can improve a little bit the upper bound for $h_G$ for some graphs. However, the upper bound is still exponential in k while we think $h_G$ is linear.

**Conjecture 27** For any graph G, there is a constant $C_G$ such that $h_G(k) \leq C_G \times k$ for any k.
A very first step towards this conjecture would be to prove it for $K_2$. Balogh et al.[1] showed that $h K_2(G) \geq 2k$ because $\Gamma(K_2 \square K_k[S_2]) = 2\Gamma(K_k[S_2]) = 2k$ with $S_2$ the edgeless graph on two vertices. They also conjectured that $h K_2(G) \geq 2k$.

Denoting by $S_k$ be the edgeless graph on k vertices, we now generalise both their conjecture and their tightness examples.

**Conjecture 28** Let k and n be two positive integers. Then $h K_n(k) = n \times k$.
More generally, we conjecture the following:

**Conjecture 29** For any graphs G and H, $\Gamma(G \square H) \leq (\Delta(H)+1)\Gamma(G)$
If true these two conjectures would be tight as shown by the following proposition.

**Proposition 30** Let k and n be two positive integers. Then $\Gamma(K_n \square K_k[S_n]) = n \times k$.

**Proof**. $\Delta(K_n \square K_k[S_n])=n\times k-1$ so $\Gamma(K_n \square K_k[S_n])\leq n\times k$.

We now prove by induction on k that $\Gamma(K_n \square K_k[S_n])\geq n\times k$. The result holds trivially when k=1. Suppose now that k>1. Let us denote the vertices of $K_k[S_n]$ by $v_{ij}, 1\leq i\leq k, 1\leq j\leq n$ so that for any i the set $\{v_{ij} | 1\leq j\leq n\}$ is stable and the vertices of $K_n$ by $x_1,...,x_n$. Let $T_1= \{(x_j, v_{1j}) | 1\leq j\leq n\}$. Then $T_1$ is a $V(G)\setminus T_1$-dominating stable set. Indeed let $(x_j, v_{i1})$ be a vertex in $V(G)\setminus T_1$. Then it is adjacent to $(x_j, v_{1j})$ if $i \neq 1$ and to $(x_1,v_{11})$ if $i \neq 1$. More generally, for $1 \leq i \leq n$, the set $T_i = \{(x_j, v_{1\,i+j-1}) | 1\leq j\leq n\}$ is a $V(G)\setminus T_i$-dominating stable set. Note that $(T_1,...,T_n)$ is a partition of $\{(x_j, v_{1i}) | 1\leq j\leq n, 1\leq i\leq n\}$ and that $K_n\square K_k[S_n]$ is isomorphic to $K_n\square K_{k-1}[S_n]$. Hence applying Proposition 5, to all the $T_i$ one after another, we obtain $\Gamma(K_n\square K_k[S_n]) \geq n+\Gamma(K_n\square K_{k-1}[S_n])$. Now the induction hypotheses is yields $\Gamma(K_n \square K_k[S_n]) \geq n\times k$.

**Theorem 31** For any graph G, $h_G(k)\geq \Gamma(G)+2k-2$.

**Proof**. Set $p=\Gamma(G)$ and $n= p+2k-2$. We will prove that $\Gamma(K_k[S_n] \square G)\geq \Gamma(G)+2k-2$. Let us denote the vertices of $K_k[S_n]$ by $v_{ij}, 1\leq i\leq k, 1\leq j \leq n$ so that for any i the set $\{v_{ij} | 1\leq j \leq n\}$. According to Lemma 8, there are two vertices x and y that receive colour $\Gamma(G)$ by some greedy colouring. Observe that in $K_k[S_n] \square G$ the $G_{vj1}, 1\leq j\leq n$ are disjoint copies of G. Hence, by Lemma, we obtain a partial greedy n-colouring of $K_k[S_n] \square G$. So $\Gamma(K_k[S_n] \square G)\geq n=\Gamma(G)+2k-2$.

### 4.2. Lower bounds

As G and H are induced subgraphs of $G\square H$ then $\Gamma(G\square H) \geq \max\{\Gamma(G), \Gamma(H)\}$. Lemma 32 Let G and H be two graphs. If $\chi(H) \leq \Delta(G)$ then $\Gamma(G\square H) \geq \Gamma(H) +1$.

**Proof**. W. l.o.g. we may assume that G and H have no isolated vertices. Let v be a vertex of G of degree $\Delta(G)$ and let $u_1,..., u_{\Delta(G)}$ its neighbours. Let $S_1,...,S_{\chi(H)}$ be the stable set of a colouring of H with $\chi(H)$ colours. The set $\cup_{i=1}^{\chi(H)} \{u_i\} \times S_i$ is a $V(H_v)$-dominating stable set. So by Proposition 5, $\Gamma(G\square H) \geq \Gamma(H) +1$.

**Corollary 33** Let G and H be two connected graphs such that $\Gamma(G) =\Gamma(H) =k$. Then $\Gamma(G\square H) \geq k+1$ unless $G=H=K_1$ or $G=H=K_2$.

**Proof**. If $\chi(H) \leq \Delta(G)$ or $\chi(G)\leq\Delta(H)$, we have the result by Lemma 32. So we may assume that $(H)= \chi(G)=\Delta(G)+1=\Delta(H)+1$. Hence by Brooks Theorem [3], G and H are complete graphs or odd cycles. If $G = H= K_k$, the result follows from Proposition 35 below.

If G and H are odd cycles, then one of $P_3\square K_2$ and $C_3\square K_2$ is an induced subgraph of $G\square H$. These graphs have grundy number 4; greedy 4-colourings are given Figure 3. So $\Gamma(G\square H) \geq 4$.

**Lemma 32** yields a direct easy proof of a result of Hoffman and Johnson [9] stating that the kdimensionnal hypercube $Q_k$ has Grundy number k+1 for k $\geq$3 and $\Gamma(Q_1) =\Gamma(Q_2) =2$. Recall that $Q_1=K_2$ and for k $\geq$2 then $Q_k=Q_{k-1}\square K_2$.

**Proposition 34** (Hoffman and Johnson [9]) for k $\geq$3, $\Gamma(Q_k) =k+1$.

**Proof**. As $\Delta(Q_k) =k$, we have $\Gamma(Q_k) \leq k+1$. Let us now prove the by induction that $\Gamma(Q_k) \geq k+1$. If k=3, a greedy 4-colouring is given. If k>3, then $\chi(K_2) \leq\Delta(Q_{k-1})=k$. Hence by Lemma33, $\Gamma(Q_k) \geq \Gamma(Q_{k-1}) +1\geq k+1$

**Proposition 35** For any p $\geq$2 then $\Gamma(K_p\square K_p) =2p-2$.

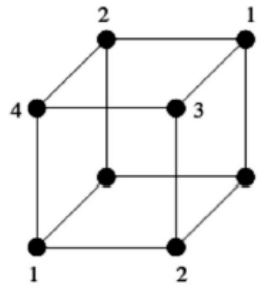

Figure 3: Partial greedy 4-colouring of the 3-dimensionnal hypercube Q

Corollary 33 implies that $\phi_\square(k,k) \geq k+1$. To have better lower bound in $\phi_\square(k,k)$, one may study the function $g(k) = \min\{\Gamma(G\square G)|\Gamma(G)=k\}$. Clearly, $g(2)=2$ and by Corollary 33 and Proposition 35, if $k\geq 3$, we have:

$k \leq \phi_\square(k,k) \leq g(k) \leq 2k-2$.

Moreover every graph with grundy number 3 has either a $K_3$ or a $P_4$ as induced subgraph. But $\Gamma(K_3\square K_3)=4$ by Proposition 35 and $\Gamma(P_4\square P_4)=5$ (As $\Delta(P_4\square P_4)\leq 4$ then $\Gamma(P_4\square P_4)\leq 5$ and it is easy to find a greedy 5-colouring of $\Gamma(P_4\square P_4)$. Hence $g(3)=4$.

## 5. DIRECT PRODUCT

A well known conjecture on graph colouring regards the chromatic number of the direct product of graphs.

**Conjecture36** (Hedetniemi[7])
$\chi(G\times H) = \min\{\chi(G), \chi(H)\}$.

Poljak [10] proved that the function f defined by $f(n)=\min\{\chi(G\times H)|\chi(G)=\chi(H)=n\}$ is either bounded by 9 or tends to infinity when n tends to infinity.

In this section, our aim is to find upper bounds of the grundy number of the direct product of two graphs in terms of the grundy number of these graphs. Ideally, we would like to determine the functions

$\phi_\times(k,l) = \min\{\Gamma(G\times H)|\Gamma(G)=k \text{ and } \Gamma(H)=l\}$.
$\Phi_\times(k,l) = \max\{\Gamma(G\times H)|\Gamma(G)=k \text{ and } \Gamma(H)=l\}$.

Let us first observe that if $G=G_1+G_2$ then $G\times H=(G_1\times H)+(G_2\times H)$. Hence it is sufficient to consider connected graphs. Furthermore, the direct product of a graph with $K_1$ is a graph without any edge of order $|G|$. So $\phi_\times(k,1)=\Phi_\times(k,1)=1$. In the remaining of this section, all the graphs are assumed to be connected of order at least 2. In particular, their grundy number is at least two.

### 5.1. Lower bounds

As every graph with grundy number k contains a k-atom as an induced subgraph then $\phi_\times(k,l)= \min\{\Gamma(G\times H)|G \text{ is a k-atom and } H \text{ is an l-atom}\}$. Furthermore if $k\geq k'$ and $l\geq l'$ then $\phi_\times(k,l)\geq \phi_\times(k',l')$.

**Theorem 37** Let G and H be two graphs with at least one edge. Then $\Gamma(G\times H) \geq \Gamma(G)+\Gamma(H)-2$. Hence if $k\geq 2$ and $l\geq 2$ then $\phi_\times(k,l)\geq k+l-2$.

**Proof**. Let k=Γ (G) and l=Γ (H). We prove the result by induction on k+l, the result holding trivially if k=l=2.

Suppose now that k+l >4. Without loss of generality, we may assume that k ≥l. Let $S_1,...,S_k$ be the stable sets of a greedy p-colouring of G. Set G'= G−$S_1$. Then $S_1$ is a (V (G'))-dominating stable set and Γ (G') =k−1. Now, in G×H), the set $S_1$×V (H) is V (G'×H)-dominating. Hence, by Proposition5, Γ (G×H) ≥Γ (G□×H). Now, since Γ (G') +Γ (H) =k+l−1, by induction hypothesis, Γ (G'×H) ≥k+l−3. So Γ (G×H) ≥k+l−2.

This lower bound for $\phi_{k,l}$ is attained when l=2 or k=l=3.

**Corollary 38** For any integer k ≥2
$$\phi\times (k, 2)=k. \quad (1)$$
$$\phi\times(3,3)=4. \quad (2)$$

**Proof**. (1) The maximum degree of $K_k \times K_2$ is k−1, so Γ ($K_k \times K_2$) ≤ k. So $\phi\times$ (k, 2) ≤ k. But Theorem37 yields $\phi\times$ (k, 2) ≥k.

(2) One can easily check that Γ ($P_4 \times P_4$) = Γ ($P_4 \times C_3$) = Γ ($C_3 \times C_3$) = 4.

## 5.2. Upper bounds

**Lemma 40** Let G and H be two graphs, and u and v be two vertices of G. If $N_G$ (u) = $N_G$ (v) then Γ (G×H) = Γ ((G−u) ×H).

Proposition3 and Lemma40 directly imply:

**Corollary 41** Let H be a graph such that Γ (H) =2. Then for every G, Γ (G×H) =Γ (G×$K_2$).
In particular, Φ (k, 2) =max {Γ (G×$K_2$)|Γ(G)=k} If G is bipartite then G×$K_2$=G+G. So, by Proposition6, Γ (G×$K_2$) =Γ (G). Then Proposition3 yields Φ× (2, 2) = 2. There are non-bipartite graphs G for which Γ (G×$K_2$) =Γ (G). For example, $K_3 \times K_2$ is the 6-cycle so Γ ($K_3 \times K_2$) =3=Γ ($K_3$). There are also graphs G for which Γ (G×$K_2$) = Γ (G): the jellyfish for example.

**Definition 42** The head of the jellyfish J is the vertex h. Let G be a graph. Then the jellyfished of G is the graph J(G) obtained from G by creating for each vertex v ∈V(G) a jellyfish J(v) whose head is identified with v.

**Proposition 43** Let G be a graph. Then Γ (J (G)) =Γ (G) +2.

**Proof**. Let c be a greedy colouring of J(G) with Γ(J(G)) colours. Let u be a vertex such that c (u) ≥4.
We claim that u is in V (G). Suppose not. Let v be the vertex of G such that u ∈J(v). Then u must be the vertex of J (v) of degree 3 adjacent to v. Since u has a neighbour of each colour smaller than 4 the vertex v must be assigned 3.

But then the two others neighbours of u are coloured 1, a contradiction.

Now as Γ (J−h) =2 the neighbours of u coloured 3 are in G. Therefore the colouring c' defined on S= {v∈V(G) | c(v)≥ 3} by c'(v)= c(v)−2 is a greedy colouring of G<S> with Γ(J(G))−2 colours. So Γ (J (G)) ≤Γ (G) +2.

**Lemma44** If G ×$K_2$ contains an induced binomial tree of index k then there is a graph H such that H ×$K_2$ contains an induced binomial tree of index k+3and Γ(H)=Γ(G)+2.

**Proof**. Let T be an induced $T_k$ of G $\times K_2$. Let S={v ∈ V(G) | {(v,1),(v,2)}⊂V($T_k$)}. Let G' be the graph obtained from G by blowing up each vertex v of S with a stable set of size two {$v_1,v_2$}. By Proposition 2 (3), Γ (G') =Γ (G).

Set H= J (G'). By Proposition 43, Γ (H) =Γ (G') +2= Γ (G) +2.

Let us now show that H $\times K_2$ contains a $T_{k+3}$. By construction, the subgraph of G'$\times K_2$ induced by (V ($T_k$) \(S$\times K_2$)) $\cup_v \cup$ S {($v_1$, 1), ($v_2$, 2)} is a$T_k$. Note that for every vertex v ∈ V (G') at most one of {(v, 1), (v, 2)} is in V(T'). Now every vertex of T' is the root of a $T_3$ in its associated J$\times K_2$. All these$T_3$ together with T' form an induced $T_{k+3}$ of H $\times K_2$.

**Corollary45** $\Phi_\times$ (2k+1, 2) ≥3k+1 and $\Phi_\times$ (2k, 2) ≥3k−1.

Note that Corollary 41 may not be generalised to graphs H such that Γ (H) = 3. Indeed 4= Γ($G_3 \times K_2$)=Γ($K_3 \times K_2$)=3

## 6. CONCLUSIONS

The Grundy number of a graph G, denoted by Γ(G), is the largest k such that G has a greedy kcoloring, that is a coloring with colours obtained by applying the greedy algorithm according to some ordering of the vertices of G. We study the Grundy number of the lexicographic, Cartesian and direct products of two graphs in terms of the Grundy numbers of these graphs .We give the lower and the upper bounds of the Grundy number for many graphs.

### ACKNOWLEDGMENTS

I thank Mr. Fréderic Havet for his collaboration.